\font\we=cmb10 at 14.4truept
\font\li=cmb10 at 10truept
\line {}
\vskip 2.5cm
\centerline {\we New Non-Abelian Zeta Functions for Curves over Finite
Fields}
\vskip 0.30cm
\centerline {\li Lin WENG}
\centerline {\bf Graduate School of Mathematics, Nagoya University, Chikusa-ku,
Nagoya 464-8602,
Japan}\vskip 0.45cm
In this paper, we introduce and study two new types of non-abelian zeta 
functions
for curves over finite fields, which are defined by using (moduli spaces of)
semi-stable vector bundles and non-stable bundles. A Riemann-Weil type 
hypothesis
is formulated for zeta functions associated to semi-stable bundles, which 
we think
is more canonical than the other one. All this is motivated by (and hence
explains in a certain sense) our work on non-abelian zeta functions for
number fields.
\vskip 0.45cm
\centerline{\bf \S 1. Restricted Non-Abelian Zeta Functions}
\vskip 0.30cm
Let $C$ be a regular, irreducible, projective curve of genus
$g$ defined over an  algebraically closed field, and $L$ a line bundle
over $C$. Then, by a result of Mumford, we know that  moduli space
${\cal M}_r(L)$ of semi-stable vector bundles of rank $r$ over $C$ with $L$ as
determinants is  projective ([Mu]). And, if $C$ and $L$ are defined over
the finite field ${\bf F}_q$ with $q$ elements, then by taking a finite 
base field
extension if necessary, we may assume that ${\cal M}_r(L)$ is also defined over
${\bf F}_q$. Moreover, by a result of Harder and Narasimhan [HN], the ${\bf
F}_q$-rational points of ${\cal M}_r(L)$ are exactly rank $r$ semi-stable 
vector
bundles of
$C$ defined over ${\bf F}_q$ with $L$ as determinants.

Clearly, the canonical line bundle $K_C$ of $C$ is defined over ${\bf F}_q$ as
well. Thus, for all
$n\in {\bf Z}$, we obtain the following natural isomorphisms  defined
over ${\bf F}_q$ too:
$$\matrix{ {\cal M}_r(L)&\to& {\cal M}_r(L\otimes K_C^{\otimes nr});&\qquad&
{\cal M}_r(L)&\to& {\cal M}_r(L^{\otimes -1}\otimes K_C^{\otimes nr})\cr
E&\mapsto&E\otimes K_C^{\otimes n};&\qquad&E&\mapsto& E^\vee\otimes 
K_C^{\otimes
n},\cr}$$  where  $E^\vee$ denotes the dual of $E$.

Define the so-called Harder-Narasimhan number $M_{C,r,L}$
by setting $$M_{C,r,L}:=\sum_{E\in {\cal M}_r(L)({\bf F}_q)}{1\over{\#{\rm
Aut}(E)}}.$$ Here ${\rm Aut}(E)$ denotes the automorphism group of $E$. 
(See e.g.,
[HN],  [Se].)  Easily, we have,  for all $n\in {\bf Z}$,
$$M_{C,r,L}=M_{C,r,L\otimes K_C^{\otimes nr}}=M_{C,r,L^{\otimes -1}\otimes
K_C^{\otimes nr}}.$$

The Harder-Narasimhan numbers may be calculated by counting  contributions of
non-stable  bundles, with the help of a famuous result of Siegel
on quadratic forms and Harder-Narasimhan's filtration for vector bundles.
(For details, see e.g., [HN], [DR] and [AB].) But in this section, we 
relate them
with a new type of (restricted) non-abelian zeta functions.

Denote by $h^i$ the dimension of cohomology groups $H^i$, $i=0,1$, by
$\chi(C,E):=h^0(C,E)-h^1(C,E)$ the Euler-Poincar\'e characteristic, and by
$d(E)$ the degree of $E$. Set
$${\cal M}_{r;L}({\bf F}_q):=\cup_{n\in {\bf Z}}\Big({\cal M}_r(L\otimes
K_C^{\otimes nr})({\bf F}_q)\cup {\cal M}_r(L^{\otimes -1}\otimes K_C^{\otimes
nr})({\bf F}_q)\Big).$$

\noindent
{\bf Definition.} {\it Let $C$ be a regular, geometrically irreducible, 
projective
curve defined over the finite field ${\bf F}_q$ with $q$ elements. With 
respect to
any fixed positive integer $r\in {\bf Z}_{>0}$, and any line bundle $L$ of $C$
defined over ${\bf F}_q$, define a weight $r$ and level $L$ restricted 
non-abelian
zeta function $\xi_{r,L}(s)$ of $C$ by
$$\xi_{r,L}(s):=\sum_{E\in {\cal M}_{r,L}({\bf
F}_q)}{{q^{h^0(C,E)}-1}\over {\#{\rm
Aut}(E)}} \cdot\Big(q^{\chi(C,E)}\Big)^{-s},\qquad {\rm Re}(s)>1.$$}
\noindent
Clearly, $\xi_{r,L}(s)
=\xi_{r,L^{\otimes \pm 1}\otimes K_C^{\otimes nr}}(s).$ Hence, from now on, we
assume
$0\leq d(L)\leq r(g-1)$.
 
\vskip 0.30cm
\noindent
{\bf Theorem 1.} {\it  For ${\rm Re}(s)>1$, we have
$$\eqalign{\xi_{r,L}(s)=&{1\over
2}\sum_{E\in {\cal M}_{r,L}({\bf F}_q);0\leq d(E)\leq
r(2g-2)}{{q^{h^0(C,E)}}\over {\#{\rm
Aut}(E)}}\cdot\Big(\big(q^{\chi(C,E)}\big)^{-s}+
\big(q^{\chi(C,E)}\big)^{s-1}\Big)\cr
&+\Big[\big(q^{d(L)-r(g-1)}\big)^{1-s}
\cdot {1\over
{q^{r(2g-2)(s-1)}-1}}+\big(q^{d(L)-r(g-1)}\big)^{s}{1\over{q^{-sr(2g-2)}-1}}\cr
&\qquad+\big(q^{d(L)-r(g-1)}\big)^{s-1}\cdot{1\over{q^{r(2g-2)(s-1)}-1}}
+\big(q^{d(L)-r(g-1))}\big)^{-s}{1\over {q^{-sr(2g-2)}-1}}\Big]\cdot
M_{C,r,L}.\cr}\eqno(*)$$ Hence, in particular, we have

\noindent
(a) $\xi_{r,L}(s)$ can be meromorphically extended to the
whole complex $t=q^{-s}$-plane;

\noindent
(b) the extension, denoted also by $\xi_{r,L}(s)$, has simple poles at $s=0$
and $s=1$  with the Harder-Narasimhan number $M_{C,r,L}$ as residues;

\noindent
(c) $\xi_{r,L}(s)$ satisfies the functional equation
$\qquad\xi_{r,L}(s)=\xi_{r,L}(1-s).$}
\vskip 0.30cm
\noindent
{\it Proof of the Theorem.} It suffices to prove (*). For this, we proceed as
follows.
Recall that for $E$  semi-stable, if
$d(E)> r(2g-2)$, then $h^1(C,E)=0$; while if $d(E)<0$, then
$h^0(C,E)=0$. Thus,
$$\eqalign{~&\xi_{r,L}(s)
=\sum_{E\in {\cal M}_{r,L}({\bf
F}_q);d(E)\geq 0}{{q^{h^0(C,E)}-1}\over {\#{\rm
Aut}(E)}}\cdot\Big(q^{\chi(C,E)}\Big)^{-s}\cr
=&\sum_{E\in {\cal M}_{r,L}({\bf
F}_q);d(E)\geq 0}{{q^{h^0(C,E)}}\over {\#{\rm
Aut}(E)}}\cdot\Big(q^{\chi(C,E)}\Big)^{-s}-\sum_{E\in {\cal M}_{r,L}({\bf
F}_q);d(E)\geq 0}{1 \over {\#{\rm
Aut}(E)}}\cdot \Big(q^{\chi(C,E)}\Big)^{-s}\cr
=&\sum_{E\in {\cal M}_{r,L}({\bf
F}_q);0\leq d(E)\leq r(2g-2)}{{q^{h^0(C,E)}}\over {\#{\rm
Aut}(E)}}\cdot\Big(q^{\chi(C,E)}\Big)^{-s}\cr
&+\sum_{E\in {\cal M}_{r,L}({\bf
F}_q);d(E)> r(2g-2)}{{q^{h^0(C,E)}}\over {\#{\rm
Aut}(E)}}\cdot\Big(q^{\chi(C,E)}\Big)^{-s}-\sum_{E\in {\cal M}_{r,L}({\bf
F}_q);d(E)\geq 0}{1\over {\#{\rm
Aut}(E)}}\cdot \Big(q^{\chi(C,E)}\Big)^{-s}.\cr}$$
So $\xi_{r,L}(s)=S_{r,L}^{0\leq d(E)\leq r(2g-2)}(s)+T_{r,L}(s)$ if we set
$$S_{r,L}^{0\leq d(E)\leq r(2g-2)}(s):=\sum_{E\in {\cal
M}_{r,L}({\bf F}_q);0\leq d(E)\leq
r(2g-2)}{{q^{h^0(C,E)}}\over {\#{\rm
Aut}(E)}}\cdot\Big(q^{\chi(C,E)}\Big)^{-s},$$ and
$$\eqalign{~T_{r,L}(s):=&\sum_{E\in {\cal M}_{r,L}({\bf
F}_q);d(E)> r(2g-2)}{{q^{h^0(C,E)}}\over {\#{\rm
Aut}(E)}}\cdot\Big(q^{\chi(C,E)}\Big)^{-s}-\sum_{E\in {\cal M}_{r,L}({\bf
F}_q);d(E)\geq 0}{1\over {\#{\rm
Aut}(E)}}\cdot \Big(q^{\chi(C,E)}\Big)^{-s}.\cr}$$
\vskip 0.30cm
\noindent
{\bf Lemma.} {\it The function $t^{r(g-1)}\cdot S_{r,L}^{0\leq d(E)\leq 
r(2g-2)}(s)$
is holomorphic in $t=q^{-s}$. Moreover $$S_{r,L}^{0\leq d(E)\leq 
r(2g-2)}(s)={1\over
2}\sum_{E\in {\cal M}_{r,L}({\bf F}_q);0\leq d(E)\leq
r(2g-2)}{{q^{h^0(C,E)}}\over {\#{\rm
Aut}(E)}}\cdot\Big(\big(q^{\chi(C,E)}\big)^{-s}+
\big(q^{\chi(C,E)}\big)^{s-1}\Big)$$ and hence
satisfies the functional equation
$$S_{r,L}^{0\leq d(E)\leq r(2g-2)}(s)=S_{r,L}^{0\leq d(E)\leq r(2g-2)}(1-s).$$}

\noindent
{\it Proof of the lemma.} Holomorphicity is clear, as only finitely many 
terms are
involved.  The functional equation comes from the Riemann-Roch theorem. Indeed,
$$\eqalign{S_{r,L}^{0\leq d(E)\leq r(2g-2)}(s)
=&{1\over
2}\Big(\sum_{E\in {\cal M}_{r,L}({\bf F}_q);0\leq d(E)\leq
r(2g-2)}{{q^{h^0(C,E)}}\over {\#{\rm
Aut}(E)}}\cdot\Big(q^{\chi(C,E)}\Big)^{-s}\cr
&\qquad+\sum_{E\in
{\cal M}_{r,L}({\bf F}_q);0\leq d(E)\leq
r(2g-2)}{{q^{h^0(C,E)}}\over {\#{\rm
Aut}(E)}}\cdot\Big(q^{\chi(C,E)}\Big)^{-s}\Big)\cr
=&{1\over 2}\Big(\sum_{E\in {\cal
M}_{r,L}({\bf F}_q);0\leq d(E)\leq
r(2g-2)}{{q^{h^0(C,E)}}\over {\#{\rm
Aut}(E)}}\cdot\Big(q^{\chi(E)}\Big)^{-s}\cr
&\qquad+\sum_{E^\vee\otimes
K_C\in {\cal M}_{r,L}({\bf F}_q);0\leq d(E^\vee\otimes K_C)\leq
r(2g-2)}{{q^{h^0(C,E^\vee\otimes K_C)}}\over {\#{\rm
Aut}(E^\vee\otimes K_C)}}\cdot\Big(q^{\chi(C,E^\vee\otimes
K_C)}\Big)^{-s}\Big)\cr
=&{1\over 2}\Big(\sum_{E\in {\cal
M}_{r,L}({\bf F}_q);0\leq d(E)\leq
r(2g-2)}{{q^{h^0(C,E)}}\over {\#{\rm
Aut}(E)}}\cdot\Big(q^{\chi(C,E)}\Big)^{-s}\cr
&\qquad+\sum_{E\in {\cal M}_{r,L}({\bf F}_q);0\leq d(E)\leq
r(2g-2)}{{q^{h^0(C,E^\vee\otimes K_C)}}\over {\#{\rm
Aut}(E)}}\cdot\Big(q^{\chi(C,E^\vee\otimes
K_C)}\Big)^{-s}\Big)\cr
=&{1\over 2}\Big(\sum_{E\in {\cal
M}_{r,L}({\bf F}_q);0\leq d(E)\leq
r(2g-2)}{{q^{h^0(C,E)}}\over {\#{\rm
Aut}(E)}}\cdot\Big(q^{\chi(C,E)}\Big)^{-s}\cr
&+\sum_{E\in {\cal M}_{r,L}({\bf F}_q);0\leq d(E)\leq
r(2g-2)}{{q^{h^1(C,E^\vee\otimes K_C)+d(E^\vee\otimes
K_C)-r(g-1)}}\over {\#{\rm
Aut}(E)}}\cdot\Big(q^{\chi(C,E^\vee\otimes K_C)}\Big)^{-s}\Big)\cr
=&{1\over 2}\Big(\sum_{E\in {\cal
M}_{r,L}({\bf F}_q);0\leq d(E)\leq
r(2g-2)}{{q^{h^0(C,E)}}\over {\#{\rm
Aut}(E)}}\cdot\Big(q^{\chi(C,E)}\Big)^{-s}\cr
&+\sum_{E\in {\cal M}_{r,L}({\bf F}_q);0\leq d(E)\leq
r(2g-2)}{{q^{h^0(C,E)-d(E)+r(g-1)}}\over {\#{\rm
Aut}(E)}}\cdot\Big(q^{-d(E)+r(g-1)}\Big)^{-s}\Big)\cr
=&{1\over
2}\sum_{E\in {\cal M}_{r,L}({\bf F}_q);0\leq d(E)\leq
r(2g-2)}{{q^{h^0(C,E)}}\over {\#{\rm
Aut}(E)}}\cdot\Big(\big(q^{\chi(C,E)}\big)^{-s}+
\big(q^{\chi(C,E)}\big)^{s-1}\Big).\cr}$$  This completes the proof of the 
lemma.
\vskip 0.30cm
As for $T_{r,L}(s)$, clearly, by the vanishing recalled at the
beginning and the Riemann-Roch, we get
$$\eqalign{T_{r,L}(s)=&\sum_{E\in {\cal M}_{r,L}({\bf F}_q);d(E)>
r(2g-2)}{{q^{d(E)-r(g-1)}}\over {\#{\rm
Aut}(E)}}\cdot\Big(q^{\chi(C,E)}\Big)^{-s}-\sum_{E\in {\cal M}_{r,L}({\bf
F}_q);d(E)\geq 0}{1\over {\#{\rm Aut}(E)}}\cdot\Big(q^{\chi(C,E)}\Big)^{-s}\cr
=&\sum_{E\in {\cal M}_{r,L}({\bf F}_q);d(E)>
r(2g-2)}{1\over {\#{\rm Aut}(E)}}\cdot\Big(q^{\chi(C,E)}\Big)^{1-s}
-\sum_{E\in {\cal M}_{r,L}({\bf F}_q);d(E)\geq
0}{1\over {\#{\rm Aut}(E)}}\cdot\Big(q^{\chi(C,E)}\Big)^{-s}\cr
=&\Big(\sum_{E\in {\cal M}_r(L\otimes K_C^{\otimes rn})
({\bf
F}_q);d(E)> r(2g-2)}+\sum_{E\in {\cal M}_r(L^{\otimes -1}\otimes 
K_C^{\otimes rn})({\bf
F}_q);d(E)> r(2g-2)}\Big)\cdot{1\over {\#{\rm
Aut}(E)}}\cdot\Big(q^{\chi(C,E)}\Big)^{1-s}\cr &-\Big(\sum_{E\in {\cal
M}_r(L\otimes K_C^{\otimes rn})({\bf F}_q);d(E)\geq 0}+\sum_{E\in {\cal
M}_r(L^{\otimes -1}\otimes K_C^{\otimes rn})({\bf F}_q);d(E)\geq
0}\Big)\cdot{1\over {\#{\rm Aut}(E)}}\cdot\Big(q^{\chi(C,E)}\Big)^{-s}\cr }$$
$$\eqalign{
=&\Big(\sum_{E\in {\cal
M}_r(L\otimes K_C^{\otimes rn}) ({\bf
F}_q);d(E)> r(2g-2)}{1\over {\#{\rm 
Aut}(E)}}\cdot\Big(q^{\chi(C,E)}\Big)^{1-s}\cr
&\qquad-\sum_{E\in {\cal
M}_r(L^{\otimes -1}\otimes K_C^{\otimes rn})({\bf F}_q);d(E)\geq
0}{1\over {\#{\rm Aut}(E)}}\cdot\Big(q^{\chi(C,E)}\Big)^{-s}\Big)\cr &+\Big(
\sum_{E\in {\cal M}_r(L^{\otimes -1}\otimes K_C^{\otimes rn})({\bf F}_q);d(E)>
r(2g-2)}{1\over {\#{\rm Aut}(E)}}\cdot\Big(q^{\chi(C,E)}\Big)^{1-s}\cr
&\qquad-\sum_{E\in {\cal M}_r(L\otimes
K_C^{\otimes rn})({\bf F}_q);d(E)\geq
0}{1\over {\#{\rm Aut}(E)}}\cdot\Big(q^{\chi(C,E)}\Big)^{-s}\Big)\cr
=&\Big[\Big(\sum_{n;d(E)=d(L)+nr(2g-2)>
r(2g-2)}\Big(q^{\chi(C,E)}\Big)^{1-s}-\sum_{n;d(E)=-d(L)+nr(2g-2)\geq
0}\Big(q^{\chi(C,E)}\Big)^{-s}\Big)\cr &\qquad+\Big(
\sum_{n;d(E)=-d(L)+nr(2g-2)>
r(2g-2)}\Big(q^{\chi(C,E)}\Big)^{1-s}-\sum_{n;d(E)=d(L)+nr(2g-2)\geq
0}\Big(q^{\chi(C,E)}\Big)^{-s}\Big)\Big]\cdot M_{C,r,L},\cr}$$
in which  we omit stating  that $E$'s are semi-stable
with determinants
$L^{\otimes \pm}\otimes K_C^{\otimes nr}$.
Therefore,
$$\eqalign{T_{r,L}(s)
=&\Big[\Big(\sum_{n=1}^\infty\Big(q^{d(L)+nr(2g-2)-r(g-1)}\Big)^{1-s}
-\sum_{n=1}^\infty\Big(q^{-d(L)+nr(2g-2)-r(g-1)}\Big)^{-s}\Big)\cr
&+\Big(\sum_{n=2}^\infty\Big(q^{-d(L)+nr(2g-2)-r(g-1)}\Big)^{1-s}
-\sum_{n=0}^\infty\Big(q^{d(L)+nr(2g-2)-r(g-1)}\Big)^{-s}\Big)\Big]\cdot
M_{C,r,L}\cr 
=&\Big[\Big(\sum_{n=1}^\infty\Big(q^{\chi(E_0)+nr(2g-2)}\Big)^{1-s}
-\sum_{n=1}^\infty\Big(q^{-\chi(C,E_0)+nr(2g-2)}\Big)^{-s}\Big)\cr
&+\Big(\sum_{n=2}^\infty\Big(q^{-\chi(C,E_0)+(n-1)r(2g-2)}\Big)^{1-s}
-\sum_{n=0}^\infty\Big(q^{\chi(C,E_0)+nr(2g-2)}\Big)^{-s}\Big)\Big]\cdot
M_{C,r,L},\cr}$$ for any $E_0\in {\cal M}_r(L)({\bf F}_q)$.    Thus, if
$q^{r(2g-2)(1-s)}< 1$ and $q^{-sr(2g-2)}<1$, by $0\leq d(L)\leq r(g-1)$,
$$\eqalign{T_{r,L}(s)
=&\Big[\Big(\sum_{n=1}^\infty\Big(q^{\chi(C,E_0)+nr(2g-2)}\Big)^{1-s}
-\sum_{n=1}^\infty\Big(q^{-\chi(C,E_0)+(n-1)r(2g-2)}\Big)^{-s}\Big)\cr
&+\Big(\sum_{n=2}^\infty\Big(q^{-\chi(C,E_0)+(n-1)r(2g-2)}\Big)^{1-s}
-\sum_{n=0}^\infty\Big(q^{\chi(C,E_0)+nr(2g-2)}\Big)^{-s}\Big)\Big]\cdot
M_{C,r,L}\cr
=&\Big[\Big(q^{\chi(C,E_0)}\Big)^{1-s}\sum_{n=1}^\infty\Big(q^{r(2g-2)(1-s)} 
\Big)^{n}
-\Big(q^{\chi(C,E_0)}\Big)^{s}\sum_{n=0}^\infty\Big(q^{-sr(2g-2)}\Big)^{n}\cr
&+\Big(q^{\chi(C,E_0)}\Big)^{s-1}\sum_{n=1}^\infty\Big(q^{r(2g-2)(1-s)}\Big) 
^{n}
-\Big(q^{\chi(C,E_0)}\Big)^{-s}\sum_{n=0}^\infty\Big(q^{-sr(2g-2)}\Big)^{n}\ 
Big]
\cdot M_{C,r,L}\cr
=&\Big[\Big(q^{\chi(C,E_0)}\Big)^{1-s}
\cdot {{q^{r(2g-2)(1-s)}}\over
{1-q^{r(2g-2)(1-s)}}}+\Big(q^{\chi(C,E_0)}\Big)^{s}{1\over{q^{-sr(2g-2)}-1}}\cr
&+\Big(q^{\chi(C,E_0)}\Big)^{s-1}\cdot{{q^{r(2g-2)(1-s)}}\over{1-q^{r(2g-2)( 
1-s)}}}
+\Big(q^{\chi(C,E_0)}\Big)^{-s}{1\over {q^{-sr(2g-2)}-1}}\Big]\cdot 
M_{C,r,L}\cr
=&\Big[\big(q^{\chi(C,E_0)}\big)^{1-s}
\cdot {1\over
{q^{r(2g-2)(s-1)}-1}}+\big(q^{\chi(C,E_0)}\big)^{s}{1\over{q^{-sr(2g-2)}-1}}\cr
&+\big(q^{\chi(C,E_0)}\big)^{s-1}\cdot{1\over{q^{r(2g-2)(s-1)}-1}}
+\big(q^{\chi(C,E_0)}\big)^{-s}{1\over {q^{-sr(2g-2)}-1}}\Big]\cdot
M_{C,r,L}.\cr}$$ This then proves the existence of meromorphic extension, the
statement for simple poles and their residues, and the functional equation for
$T_{r,L}(s)$. Thus by the Lemma, we complete the proof of the Theorem.
\vskip 0.45cm
\centerline {\bf \S 2. New Non-Abelian Zeta Functions for Curves over 
Finite Fields}
\vskip 0.30cm
 From definition, clearly,  $${\xi_{1,L}(s)=q^{(g-1)s}\sum_{E=L^{\otimes
\pm}\otimes K_C^{\otimes n};n\in {\bf
Z}}{{q^{h^0(C,E)}-1}\over {q-1}}\cdot(q^{d(E)})^{-s}.}$$ So, by taking a 
finite sum
over suitable  $L$, we could arrive at
$${q^{(g-1)s}\sum_{M\in {\rm Pic}(C)({\bf
F}_q)}{{q^{h^0(C,E)}-1}\over {q-1}}\cdot(q^{d(E)})^{-s},}$$ which is 
nothing but
$$q^{(g-1)s}\cdot \sum_{D\geq
0}N(D)^{-s},$$ i.e.,  the {\it standard abelian zeta function} $\zeta_C(s)$, or
Artin zeta function ([A]), for $C$ times with $q^{(g-1)s}$.
This then suggests the the following discussion.

Denote by ${\cal M}_{r,d}(C)$ the moduli space of degree $d$ semi-stable
vector bundles of rank $r$ on $C$. Fixed a degree 1 line bundle $A$ on
$C$ which is ${\bf F}_q$-rational. (One may give a proof of this fact by using
properties of Artin zeta function $\zeta_C(s)$.)
Clearly, we then have the following isomorphisms defined over ${\bf F}_q$:
$$\matrix{ {\cal M}_r(L)&\to& {\cal M}_r(L\otimes A^{\otimes
nr})\cr E&\mapsto&E\otimes A^{\otimes
n}.\cr}$$ Thus, as before,  we may assume that, if necessary, by taking a 
finite
fields extension, ${\cal M}_{r,d}(C)$ are defined over ${\bf F}_q$
as well. Moreover, we know that all ${\bf F}_q$-rational semi-stable
vector bundles are indeed an element in
$${\cal M}_r(C)({\bf F}_q):=\cup_d{\cal M}_{r,d}(C)({\bf F}_q).$$
\noindent
{\bf Definition.} {\it For $r\in {\bf Z}_{>0}$, define a weight $r$
non-abelain zeta function $\zeta_{C,r}(s)$ of
$C$ by
$$\zeta_{C,r}(s):=\sum_{E\in {\cal M}_{r}(C)({\bf
F}_q)}{{q^{h^0(C,E)}-1}\over{\#{\rm
Aut}(E)}}\cdot\Big(q^{d(E)}\Big)^{-s},\qquad {\rm Re}(s)>1.$$}
Set also $\xi_{C,r}(s):=q^{r(g-1)s}\cdot\zeta_{C,r}(s)$, $t=q^{-s}$ and
$Z_{C,r}(t)=\zeta_{C,r}(s)$.
\vskip 0.30cm
\noindent
{\bf Theorem 2.} {\it (a) $\zeta_{C,r}(s)$ is absolutely convergent for
${\rm Re}(s)>1$, can be meromorphically extended to the whole complex 
$t=q^{-s}$
plane which has  simple poles at $s=0$ and
$s=1$  with the same residue, and satisfies the
functional equation
$\qquad\xi_{r,L}(s)=\xi_{r,L}(1-s).$

\noindent
(b) $\displaystyle{Z_{r,C}(t)={{P(t)}\over{(1-t^r)(1-(qt)^r)}}}$ where 
$P(t)\in {\bf
Q}[t]$ is a degree $2rg$ polynomial with rational coefficients.

\noindent
(c) Let $\displaystyle{P(t)=P(0)\cdot \prod_{i=1}^{2rg}(1-\omega_it),}$
then after a suitable arrangement, we have
$$\omega_i\cdot\omega_{2rg-i}=q,\qquad i=1,\dots, rg.$$}
\noindent
 From the proof  given below, one can give the precise values of the
residues in terms of a certain combination of Harder-Narasimhan numbers. We 
leave
this to the reader. Moreover, (c) suggests that  Weil-Riemann Hypothesis ([W])
holds also for our new zeta functions. This then leads to the following
\vskip 0.30cm
\noindent
{\bf Riemann-Weil Hypothesis.} {\it The reciprocal roots $\omega_i$,
$i=1,\dots,2rg$, of the weight
$r$ non-abelain zeta functions of curves defined over finite
fields satisfy $$|\omega_i|=q^{1\over 2},\qquad i=1,\dots, 2rg.$$}
\vskip 0.30cm
\noindent
{\it Proof of the Theorem.} By a finite field extension, we may assume that
all ${\cal M}_r(L)$ are defined over ${\bf F}_q$ if $L$ is defined over 
${\bf F}_q$.
With this,  (a) is a direct consequence of Theorem  1 in \S1 by the fact that
$${\cal M}_{r,d}(C)({\bf F}_q)=\cup_{d=0}^{r-1}\cup_{L\in {\rm 
Pic}^d(C)({\bf F}_q)}
\cup_{k=1}^{g-1}{\cal M}_{r,L\otimes A^{kr}}({\bf F}_q).$$ Hence we only 
need to
prove (b) and (c).

Let us assume (b), then by the functional equation for $t=q^{-s}$, which 
changes
$t$ to ${1\over {qt}}$, $\prod_{i=1}^{2rg}(1-\omega_it)$ changes to
$\prod_{i=1}^{2rg}(1-{{\omega_i}\over {qt}})$ accordingly.
Thus
$$\prod_{i=1}^{2rg}(1-\omega_it)=\prod_{k=1}^{2rg}(1-{1\over
{\omega_kq}}t).$$ This then implies (c).

So it suffices to prove (b), the rationality of our non-abelian zeta function.
For this, let us first assume that
$$Z_{r,C}(t)={{P(t)}\over{(1-t^r)(1-(qt)^r)}}$$ where $P(t)\in {\bf
Q}[t]$. Then by the functional equation for $\zeta_{r,C}(s)$, we have
$${{P(t)}\over{(1-t^r)(1-(qt)^r)}}\cdot t^{-r(g-1)}=
{{P({1\over {qt}})}\over{(1-{1\over {(qt)^r}})(1-{1\over {t^r}})}}\cdot
(qt)^{r(g-1)}.$$
Thus, ${\rm deg}(P)-2r-r(g-1)=r(g-1)$ by comparing degrees of $t$ of
rational functions  on both sides. That is to say, we have ${\rm
deg}(P)=2rg$. In this way, we are lead to prove the following
\vskip 0.30cm
\noindent
{\bf Theorem$'$.} {\it With the same notation as above, there exists a 
polynomial
$P(t)\in {\bf Q}[t]$ such that
$$Z_{r,C}(t)={{P(t)}\over{(1-t^r)(1-(qt)^r)}}.$$}

\noindent
{\it Proof.} Note that we have the following vanishing:
for $E$  semi-stable, if
$d(E)> r(2g-2)$, then $h^1(C,E)=0$; while if $d(E)<0$, then
$h^0(C,E)=0$. Thus, by definition,
$$\eqalign{~&Z_{r,C}(t)\cr
=&\sum_{E\in {\cal M}_{r}(C)({\bf
F}_q),d(E)\geq 0}{{q^{h^0(C,E)}-1}\over{\#{\rm
Aut}(E)}}\cdot\Big(q^{d(E)}\Big)^{-s}\cr
=&\sum_{E\in {\cal M}_{r}(C)({\bf
F}_q),d(E)\geq 0}{{q^{h^0(C,E)}}\over{\#{\rm
Aut}(E)}}\cdot\Big(q^{d(E)}\Big)^{-s}-\sum_{E\in {\cal M}_{r}(C)({\bf
F}_q),d(E)\geq 0}{{1}\over{\#{\rm
Aut}(E)}}\cdot\Big(q^{d(E)}\Big)^{-s}\cr
=&\sum_{E\in {\cal M}_{r}(C)({\bf
F}_q),r(2g-2)\geq d(E)\geq 0}{{q^{h^0(C,E)}}\over{\#{\rm
Aut}(E)}}\cdot\Big(q^{d(E)}\Big)^{-s}+\sum_{E\in {\cal M}_{r}(C)({\bf
F}_q),d(E)>r(2g-2)}{{q^{h^0(C,E)}}\over{\#{\rm
Aut}(E)}}\cdot\Big(q^{d(E)}\Big)^{-s}\cr
&\qquad-\sum_{E\in {\cal M}_{r}(C)({\bf
F}_q),d(E)\geq 0}{{1}\over{\#{\rm
Aut}(E)}}\cdot\Big(q^{d(E)}\Big)^{-s}.\cr}$$
Set $$\eqalign{I:=&\sum_{E\in {\cal M}_{r}(C)({\bf
F}_q),r(2g-2)\geq d(E)\geq 0}{{q^{h^0(C,E)}}\over{\#{\rm
Aut}(E)}}\cdot\Big(q^{d(E)}\Big)^{-s}\cr
II:=&\sum_{E\in {\cal M}_{r}(C)({\bf
F}_q),d(E)>r(2g-2)}{{q^{h^0(C,E)}}\over{\#{\rm
Aut}(E)}}\cdot\Big(q^{d(E)}\Big)^{-s}\cr
III:=&\sum_{E\in {\cal M}_{r}(C)({\bf
F}_q),d(E)\geq 0}{{1}\over{\#{\rm
Aut}(E)}}\cdot\Big(q^{d(E)}\Big)^{-s}\cr}$$ Then
$$Z_{r,C}(t)=I(t)+II(t)+III(t).$$ Thus it suffices to prove the following
\vskip 0.30cm
\noindent
{\bf Lemma.} {\it With the same notation as above, $$I(t),\quad
\Big(1-(qt)^r\Big)\cdot II(t),\quad {\rm and}\quad \Big(1-t^r\Big)\cdot
III(t)\in {\bf Q}[t].$$}
\noindent
{\it Proof of the Lemma.} By definition, $I(t)\in {\bf Q}[t]$.
Hence we should prove the assertions for $II(t)$ and $III(t)$.

We first deal with $II(t)$. By the vanishing, we have $h^0(C,E)=d-r(g-1)$. 
Hence
$$II(t)=q^{-r(g-1)}\sum_{E\in {\cal M}_{r}(C)({\bf
F}_q),d(E)>r(2g-2)}{{(qt)^{d(E)}}\over{\#{\rm
Aut}(E)}}.$$ Let ${\cal M}_{r,k}:={\cal M}_{r,k}(C)({\bf F}_q)$,
then by tensoring with $A^\otimes n$, $n\in {\bf Z}_{\geq 0}$, we have
isomorphisms
$$\eqalign{~&{\cal M}_{r,0}\simeq {\cal M}_{r,r}\simeq\dots\simeq{\cal
M}_{r,r(2g-2)}\simeq {\cal
M}_{r,(r+1)(2g-2)}\simeq\dots\cr
&{\cal M}_{r,1}\simeq {\cal
M}_{r,r+1}\simeq\dots\simeq{\cal M}_{r,r(2g-2)+1}\simeq {\cal
M}_{r,(r+1)(2g-2)+1}\simeq\dots\cr
&\dots\qquad\dots\qquad\dots\cr
&{\cal M}_{r,r-1}\simeq {\cal M}_{r,r+(r-1)}\simeq\dots\simeq{\cal
M}_{r,r(2g-2)+(r-1)}\simeq {\cal
M}_{r,(r+1)(2g-2)+(r-1)}\simeq...\cr}$$
Thus, if we further set
$$M_r:={\cal M}_{r,r(2g-2)+1}\cup\dots {\cal
M}_{r,(r+1)(2g-2)+(r-1)}\cup {\cal
M}_{r,(r+1)(2g-2)},$$
then
$$\eqalign{II(t)
=&q^{-r(g-1)}\sum_{E\in
M_r}\sum_{n=0}^\infty{{(qt)^{d(E\otimes A^{\otimes n})}}\over{\#{\rm
Aut}(E)}}\cr
=&q^{-r(g-1)}\sum_{E\in
M_r}\sum_{n=0}^\infty{{(qt)^{d(E)+nr}}\over{\#{\rm
Aut}(E)}}=q^{-r(g-1)}\sum_{E\in
M_r}{{(qt)^{d(E)}}\over{\#{\rm
Aut}(E)}}\cdot\sum_{n=0}^\infty((qt)^r)^n
\cr
=&q^{-r(g-1)}\sum_{E\in
M_r}{{(qt)^{d(E)}}\over{\#{\rm
Aut}(E)}}\cdot{1\over{1-(qt)^r}}.\cr}$$ This proves the assertion for
$II(t)$.

Similarly, for $III(t)$, set $$M_r':= {\cal M}_{r,0}\cup {\cal
M}_{r,1}\cup\dots\cup {\cal M}_{r,r-1}.$$ Then
$$\eqalign{III(t)
=&\sum_{E\in
M_r'}\sum_{n=0}^\infty{{t^{d(E\otimes A^{\otimes n})}}\over{\#{\rm
Aut}(E)}}\cr
=&\sum_{E\in
M_r'}\sum_{n=0}^\infty{{t^{d(E)+nr}}\over{\#{\rm
Aut}(E)}}=\sum_{E\in
M_r}{{t^{d(E)}}\over{\#{\rm
Aut}(E)}}\cdot\sum_{n=0}^\infty(t^r)^n\cr
=&\sum_{E\in
M_r'}{{t^{d(E)}}\over{\#{\rm
Aut}(E)}}\cdot{1\over{1-t^r}}.\cr}$$ This proves the assertion for
$III(t)$. All in all, we have proved the lemma, hence the theorem$'$ and the
theorem itself.
\vskip 0.30cm
As a direct consequence of the rationality of our new zeta functions, i.e., 
(b) of
the Theorem, we have the following

\noindent
{\bf Corollary.} {\it (a) For each $m\geq 1$ there exists suitable number 
$N_m$ such
that
$$Z_{r,C}(t)=\exp\Big(\sum_{m=1}^\infty N_m{{t^m}\over n}\Big).\eqno(**)$$
Moreover, $$N_m=\cases{r(1+q^m)-\sum_{i=1}^{2rg}\omega_i^m,& if $r\ |m$;\cr
-\sum_{i=1}^{2rg}\omega_i^m,& if $r\not| m$.\cr}$$

\noindent
(b) For any positive integer $a$ such that $(a,r)=1$, if $\zeta_{i,a}, 
i=1,\dots, a$
denote all $a$-th roots of unity, then $$\prod_{i=1}^aZ_{C,r}(\zeta_{i,a}t)=
\exp\Big(\sum_{m=1}^\infty N_{ma}{{T^m}\over m}\Big)\eqno(***)$$ where 
$T=t^a$.}

\noindent
{\it Proof.} In fact, note that $\log(1-x)=-\sum_{m=1}^\infty{{x^m}\over 
m}$, from
the rationality of $Z_{C,r}(t)$, we get (**) directly. As for the precise 
formula
of $N_m$, we use the following fact which also implies (b) directly:
$$\sum_{i=1}^a\zeta_{i,a}^m=\cases{a,& if $a\ |m$,\cr
0, & if $a\not|m$.\cr}$$

Clearly, (***) suggests that 
$$\prod_{i=1}^aZ_{C,r}(\zeta_{i,a}t)=Z_{C_{ar},r}(T)$$
with
$C_{ar}$ is obtained from $C$ by taking simply extension of constant fields 
from
${\bf F}_q$ to ${\bf F}_{q^{ar}}$. (Question: Why $ar$ not simply $a$?) As 
a matter
of fact,
$N_1,\dots,N_r$ do have Diophantine interpretations. Moreover, there are many
important relations in this enumerative aspect of moduli spaces of semi-stable
vector bundles. As all are associated to the above  Riemann-Weil 
hypothesis, which
we cannot verify now,   we leave them to some other occasions.
\vskip 0.45cm
\centerline {\bf \S 3. Non-Stable Contributions}
\vskip 0.30cm
In this section,  we justify our new non-abelian zeta functions
by looking at non-stable contributions. But for simplicity, we only study 
the case
when $r=2$.

So let $W_d(C)$ be the collection of isomorphism classes of rank two degree
$d$ non-stable vector bundles $E$ (over $C$) defined over ${\bf F}_q$. Thus 
clearly,
we have then  bijections  $W_0(C)\simeq W_{2n}(C)$ and $W_1(C)(\simeq 
W_{2n+1}(C)$
for all $n\in {\bf Z}$ as for semi-stable vector bundles.

Naturally, one may try to define a
new zeta function for $C$ by considering a formal summation
$$\sum_{E\in W_d(C),d\in {\bf Z}}{{q^{h^0(C,E)}-1}\over {\#{\rm
Aut}(E)}}\cdot\Big(q^{d(E)}\Big)^{-s},\qquad{\rm Re}(s)>1.$$
But this time, we have some troubles. Recall that in semi-stable case,
despite that the formal summation is for all $d\in {\bf Z}$, finally, only 
terms
with $d\geq 0$ contribute, since $q^{h^0(C,E)}-1=0$ when $d<0$. Clearly, for
non-stable bundles, this does not hold. Hence, we should modify our 
definition as
follows.

\noindent
{\bf Definition.} {\it Let $C$ be a regular, geometrically irreducible, 
projective
curve defined over the finite field ${\bf F}_q$ with $q$ elements.
Define a new zeta functions $\zeta_{C,{\rm ns}}(s)$ by setting
$$\zeta_{C,{\rm ns}}(s):=\sum_{E\in W_d(C),d\geq 0}{{q^{h^0(C,E)}-1}\over 
{\#{\rm
Aut}(E)}}\cdot\Big(q^{d(E)}\Big)^{-s},\qquad{\rm Re}(s)>1.$$}

The main result of this section will be the following

\noindent
{\bf Theorem 3.} {\it $\zeta_{C,{\rm ns}}(s)$ is well-defined and admits a
meromorphic extension to the whole complex $t=q^{-s}$-plane,  which is indeed
rational too, i.e., can be written as a quotient of two polynomials with 
rational
coefficients.}

\noindent
{\it Proof.} Call any element $E\in W_{0}(C)$ an $E_0$ and 
$E_{2n}=E_0\otimes A^n$
with $A$ as in the previous section. Similarly, call $E\in W_{1}(C)$ an 
$E_1$ and
$E_{2n+1}=E_1\otimes A^n$. Then clearly $\#{\rm Aut}(E_*)=\#{\rm 
Aut}(E_{2n+*})$.
Moreover, for $E_*\in W_*(C)$, $*=0,1$, by using Harder-Narasimhan 
filtration, over
${\bf F}_q$, we have the following short exact sequences:
$$0\to L_{+,*}\to E_*\to L_{-,*}\to 0\eqno(*)$$ such that

\noindent
($a_0$) if $E_0\in W_0(C)$, $d_{+,0}:=d(L_{+,0})=1,2,3,\dots$ and
$d_{-,0}:=d(L_{-,0})=-d(L_{+,0})=-1,-2,-3,\dots$;

\noindent
($a_1$) if $E_1\in W_1(C)$, $d_{+,1}:=d(L_{+,1})=1,2,3,\dots$ and
$d_{-,1}:=d(L_{-,1})=1-d(L_{+,0})=0,-1,-2,\dots$.

Thus if we set $L_{+,2n+*}=L_{+,*}\otimes A^n$ with $*=0,1$, then we have also

\noindent
($a_{2n}$) if $E_{2n}\in W_{2n}(C)$, 
$d_{+,2n}:=d(L_{+,2n})=n+1,n+2,n+3,\dots$ and
$d_{-,2n}:=d(L_{-,2n})=n-1,n-2,n-3,\dots$;

\noindent
($a_{2n+1}$) if $E_{2n+1}\in W_{2n+1}(C)$,
$d_{+,2n+1}:=d(L_{+,2n+1})=n+1,n+2,n+3,\dots$ and
$d_{-,2n+1}:=d(L_{-,2n+1})=n,n-1,n-2,\dots$.

Next, we compute $\#{\rm Aut}(E)$ according to whether (*) is trivial or not.

\noindent
(1) If $E_m=L_{+,m}\oplus L_{-,m}$, the automorphisms consist of ${\bf
F}_q^*\times {\bf F}_q^*$ together with the unipotents of the form $1+\phi$ 
with
$$\phi\in {\rm Hom}(L_{-,m}, L_{+,m})=H^0(C,L_{-,m}^{\otimes
-1}\otimes L_{+,m})=\cases{H^0(C,L_{-,0}^{\otimes
-1}\otimes L_{+,0}),& if $m=2n$;\cr
H^0(C,L_{-,1}^{\otimes
-1}\otimes L_{+,1}),& if $m=2n+1$.\cr}$$
Hence $$\#{\rm Aut}(E_m)=(q-1)^2\cdot \cases{q^{h^0(C,L_{-,0}^{\otimes
-1}\otimes L_{+,0})},& if $m=2n$;\cr
q^{h^0(C,L_{-,1}^{\otimes
-1}\otimes L_{+,1})},& if $m=2n+1$.\cr}$$

\noindent
(2) For non-trivial extensions, we have only one copy of ${\bf F}_q^*$ and 
hence
$$\#{\rm Aut}(E_m)=(q-1)\cdot \cases{q^{h^0(C,L_{-,0}^{\otimes
-1}\otimes L_{+,0})},& if $m=2n$;\cr
q^{h^0(C,L_{-,1}^{\otimes
-1}\otimes L_{+,1})},& if $m=2n+1$.\cr}$$ But non-trivial extensions $E_m$
correspond to non-zero elements of $H^1(C,L_{-,m}^{\otimes
-1}\otimes L_{+,m})$ and proportional vectors give isomorphic bundles. 
Hence the
number of isomorphism classes of bundles $E_m$ for which $E$ is non-trivial is
$${{q^{h^1(C, C,L_{-,m}^{\otimes -1}\otimes L_{+,m})}-1}\over{q-1}}.$$

Thus in particular, in the summation
$$\sum_{E\in W_d(C),d\geq 0}{1\over{\#{\rm Aut}(E)}}\Big(q^{d(E)}\Big)^{-s},$$
the contributions arising from  given $L_{+,2n+*}, *=0,1$  are
$$\eqalign{~&\Big({1\over {(q-1)^2q^{h^0(C,L_{-,*}^{\otimes
-1}\otimes L_{+,*})}}}+
{{q^{h^1(C,L_{-,*}^{\otimes
-1}\otimes L_{+,*})}-1}\over {(q-1)^2q^{h^0(C,L_{-,*}^{\otimes
-1}\otimes L_{+,*})}}}\Big)\cdot (q^{-s})^{2n+*}\cr
=&{{q^{h^1(C,L_{-,*}^{\otimes
-1}\otimes L_{+,*})}}\over {(q-1)^2q^{h^0(C,L_{-,*}^{\otimes
-1}\otimes L_{+,*})}}}\cdot (q^{-s})^{2n+*}\cr
=&q^{-d_{+,*}+d_{-,*}+(g-1)}\cdot (q^{-s})^{2n+*}\cr
=&\cases{q^{-2d_{+,0}+(g-1)}\cdot (q^{-s})^{2n},& if $m=2n$;\cr
q^{-2d_{+,0}+g}\cdot (q^{-s})^{2n+1},& if $m=2n+1$.\cr}\cr}\eqno(*)$$
(This trick is first used by Atiyah and Bott in [AB] p.595.)

Therefore, if $J_0(C)$ denotes degree zero Jacobian of $C$, then
$$\eqalign{~&\sum_{E\in W_d(C),d\geq 0}{1\over{\#{\rm
Aut}(E)}}\Big(q^{d(E)}\Big)^{-s}\cr
=&\#J_0(C)({\bf
F}_q)\sum_{d_{+,*}=1}^\infty\sum_{n\geq 0}^\infty
\cdot{1\over {(q-1)^2}}\cdot \Big(q^{-2d_{+,0}+(g-1)}\cdot (q^{-s})^{2n}+
q^{-2d_{+,0}+g}\cdot (q^{-s})^{2n+1}\Big)\cr
=&\#J_0(C)({\bf F}_q)\cdot{1\over {(q-1)^2}}\cdot 
{{q^{-2}}\over{1-q^{-2}}}\cdot
q^{g-1}\Big(1+q^{-2s}+q^{-4s}+\dots+q\big(
q^{-s}+q^{-3s}+q^{-5s}+\dots\big)\Big)\cr
=&{{q^{g-1}}\over{(q-1)^2\cdot(q^2-1)}}\cdot {1\over
{1-q^{-2s}}}\cdot(1+q^{1-s})\cdot \#J_0(C)({\bf F}_q),\cr}\eqno(**)$$
provided that ${\rm Re}(s)>1$.

Before going further, we reminder the reader that the summation here involves a
{\it double} infinite summations, and hence is quite different in nature 
comparing
with  semi-stable cases.

Next, consider the part $$\sum_{E\in W_d(C),d\geq 0}{{q^{h^0(C,E)}}\over{\#{\rm
Aut}(E)}}\Big(q^{d(E)}\Big)^{-s}.$$
Similarly as above, we want to group non-trivial extensions and trivial 
extensions
together so that automorphisms may be calculated easily. Clearly for this 
purpose,
then we should assume that $h^0(C,E_m)=h^0(C,L_{+,m})+h^0(C,L_{-,m})$, which is
clearly the case if $h^1(C,L_{+,m})=0$, or even better, $d_{+,m}> 2g-2$. By the
complete list for degrees at the beginning of this proof, we see that there are
only finitely many $E\in W_d(C), d\geq 0$ such that $d_{+,m}\leq 2g-2$. 
Hence, by
definition, up to a polynomial, we are essentially dealing with non-stable 
bundles
such that $h^0(C,E_m)=h^0(C,L_{+,m})+h^0(C,L_{-,m})$.
Moreover, if we can calculate $h^0(C,L_{+,m})$ and $h^0(C,L_{-,m})$ 
precisely in
terms of degrees $d_{+,m}$ and $d_{-,m}$, say in  cases when we have 
vanishing of
$h^0$ or $h^1$, proceeding similarly as in (**), we can prove the theorem.

With this in mind, we write
$$\eqalign{&\sum_{E\in W_d(C),d\geq 0}{{q^{h^0(C,E)}}\over{\#{\rm
Aut}(E)}}\Big(q^{d(E)}\Big)^{-s}\cr
=&\sum_{E\in E_*\otimes A^n,E_*\in
W_*(C),*=0,1,n\leq 2g-1}{{q^{h^0(C,E)}}\over{\#{\rm
Aut}(E)}}\Big(q^{d(E)}\Big)^{-s}\cr
&\qquad+\sum_{E\in E_*\otimes A^n,E_*\in
W_*(C),*=0,1,n\geq 2g}{{q^{h^0(C,L_{+,2n+*})+h^0(C,L_{-,2n+*})}}\over{\#{\rm
Aut}(E)}}\Big(q^{d(E)}\Big)^{-s}\cr
=&\sum_{E\in E_*\otimes A^n,E_*\in
W_*(C),*=0,1,n\leq 2g-1}{{q^{h^0(C,E)}}\over{\#{\rm
Aut}(E)}}\Big(q^{d(E)}\Big)^{-s}\cr
&\qquad
+\sum_{E\in E_*\otimes A^n,E_*\in W_*(C),*=0,1,n\geq
2g,0\leq d_{-,m}\leq 
2g-2}{{q^{h^0(C,L_{+,2n+*})+h^0(C,L_{-,2n+*})}}\over{\#{\rm
Aut}(E)}}\Big(q^{d(E)}\Big)^{-s}\cr
+&\sum_{E\in E_*\otimes A^n,E_*\in W_*(C),*=0,1,n\geq
2g,d_{-,m}> 2g-2}{{q^{\chi(C,E_m)}}\over{\#{\rm
Aut}(E)}}\Big(q^{d(E)}\Big)^{-s}\cr
&\qquad
+\sum_{E\in E_*\otimes A^n,E_*\in W_*(C),*=0,1,n\geq
2g, d_{-,m}<0}{{q^{\chi(C,L_{+,m})}}\over{\#{\rm
Aut}(E)}}\Big(q^{d(E)}\Big)^{-s},\cr}$$ by a simple observation on cohomology
groups. Set now
$$\eqalign{I:=&\sum_{E\in E_*\otimes A^n,E_*\in W_*(C),*=0,1,n\geq
2g,d_{-,m}> 2g-2}{{q^{\chi(C,E_m)}}\over{\#{\rm
Aut}(E)}}\Big(q^{d(E)}\Big)^{-s},\cr
II:=&\sum_{E\in E_*\otimes A^n,E_*\in W_*(C),*=0,1,n\geq
2g, d_{-,m}<0}{{q^{\chi(C,L_{+,m})}}\over{\#{\rm
Aut}(E)}}\Big(q^{d(E)}\Big)^{-s},\cr
III:=&\sum_{E\in E_*\otimes A^n,E_*\in W_*(C),*=0,1,n\geq
2g,0\leq d_{-,m}\leq 
2g-2}{{q^{h^0(C,L_{+,2n+*})+h^0(C,L_{-,2n+*})}}\over{\#{\rm
Aut}(E)}}\Big(q^{d(E)}\Big)^{-s},\cr
IV:=&\sum_{E\in E_*\otimes A^n,E_*\in
W_*(C),*=0,1,n\leq 2g-1}{{q^{h^0(C,E)}}\over{\#{\rm
Aut}(E)}}\Big(q^{d(E)}\Big)^{-s}.\cr}$$ Then clearly,

\noindent
(i) $\displaystyle{\sum_{E\in W_d(C),d\geq 0}{{q^{h^0(C,E)}}\over{\#{\rm
Aut}(E)}}\Big(q^{d(E)}\Big)^{-s}=I+II+III+IV.}$

\noindent
(ii) By a similar calculation as in (**) for the case
$\sum_{E}{{1}\over{\#{\rm
Aut}(E)}}\Big(q^{d(E)}\Big)^{-s}$, we see
  that $I$ and $II$ are convergent provided ${\rm Re}(s)>1$ and
are rational.

Therefore, it suffices to deal with $III$ and $IV$.  We study $IV$ first. 
For this,
note that $d_{+,m}+d_{-,m}=m$ or $m+1$,
we have then
$$\eqalign{IV=&\sum_{E\in E_*\otimes A^n,E_*\in
W_*(C),*=0,1,n\leq 2g-1,d_{+,m}\leq 2g-2}{{q^{h^0(C,E)}}\over{\#{\rm
Aut}(E)}}\Big(q^{d(E)}\Big)^{-s}\cr
&\qquad
+\sum_{E\in E_*\otimes A^n,E_*\in
W_*(C),*=0,1,n\leq 2g-1,d_{+,m}>
2g-2}{{q^{h^0(C,L_{+,m})+h^0(C,L_{-,m})}}\over{\#{\rm
Aut}(E)}}\Big(q^{d(E)}\Big)^{-s}\cr
=&\sum_{E\in E_*\otimes A^n,E_*\in
W_*(C),*=0,1,n\leq 2g-1,d_{+,m}\leq 2g-2}{{q^{h^0(C,E)}}\over{\#{\rm
Aut}(E)}}\Big(q^{d(E)}\Big)^{-s}\cr
&\qquad
+\sum_{E\in E_*\otimes A^n,E_*\in
W_*(C),*=0,1,n\leq 2g-1,d_{+,m}>
2g-2,0\leq d_{-,m}\leq 2g-1}{{q^{h^0(C,L_{+,m})+h^0(C,L_{-,m})}}\over{\#{\rm
Aut}(E)}}\Big(q^{d(E)}\Big)^{-s}\cr
&\qquad
+\sum_{E\in E_*\otimes A^n,E_*\in
W_*(C),*=0,1,n\leq 2g-1,d_{+,m}>
2g-2,d_{-,m}< 0}{{q^{\chi(C,L_{+,m})}}\over{\#{\rm
Aut}(E)}}\Big(q^{d(E)}\Big)^{-s}.\cr}$$ Now set
$$\eqalign{IV_a:=&\sum_{E\in E_*\otimes A^n,E_*\in
W_*(C),*=0,1,n\leq 2g-1,d_{+,m}\leq 2g-2}{{q^{h^0(C,E)}}\over{\#{\rm
Aut}(E)}}\Big(q^{d(E)}\Big)^{-s}\cr
IV_b:=&\sum_{E\in E_*\otimes A^n,E_*\in
W_*(C),*=0,1,n\leq 2g-1,d_{+,m}>
2g-2,0\leq d_{-,m}\leq 2g-1}{{q^{h^0(C,L_{+,m})+h^0(C,L_{-,m})}}\over{\#{\rm
Aut}(E)}}\Big(q^{d(E)}\Big)^{-s}\cr
IV_c:=&\sum_{E\in E_*\otimes A^n,E_*\in
W_*(C),*=0,1,n\leq 2g-1,d_{+,m}>
2g-2,d_{-,m}< 0}{{q^{\chi(C,L_{+,m})}}\over{\#{\rm
Aut}(E)}}\Big(q^{d(E)}\Big)^{-s}.\cr}$$
Then

\noindent
(a) $IV=IV_a+IV_b+IV_c$;

\noindent
(b) $IV_a$ and $IV_b$ consist of only finitely many terms and hence are
polynomials;

\noindent
(c) $IV_c$ may be calculated similarly as in the case $\sum_{E}{{1}\over{\#{\rm
Aut}(E)}}\Big(q^{d(E)}\Big)^{-s}$.

Hence we get the convergence provided ${\rm Re}(s)>1$ and the rationality 
for $IV$.

Finally, let us deal with $III$. This seems to be a rather delicate work: the
summation is an infinite one; while previous trick cannot be applied
since we cannot calculate $h^0(C,L_{-,m})$ precisely.  However, we are quite
lucky: In any case, now we have
$h^0(C,E_m)=h^0(C,L_{+,m})+h^0(C,L_{-,m})$ and
$h^0(C,L_{+,m})=\chi(C,L_{+,m})$. Hence,
$$III=\sum_{E\in E_*\otimes A^n,E_*\in W_*(C),*=0,1,n\geq
2g,0\leq d_{-,m}\leq
2g-2}{{q^{h^0(C,L_{-,m})+d_{-,m}}}\over{(q-1)^2}}\Big(q^{2n+*}\Big)^{-s}.$$
Set now for $*=0,1$,
$$A_*(m):=\sum_{E\in E_*\otimes A^n,E_*\in W_*(C),n\geq 2g,0\leq d_{-,m}\leq
2g-2}{{q^{h^0(C,L_{-,m})+d_{-,m}}}\over{(q-1)^2}}.$$
Then we have $$A_*(m)=\sum_{L\in {\rm Pic}^d(C)({\bf F}_q),0\leq d\leq 2g-2}
{{q^{h^0(C,L)}+d(L)}\over {(q-1)^2}}:=A(C)$$ which is independent of $*$ and
$m$. Therefore
$$III=A(C)\times \sum_{*=0,1,n\geq 2g}\Big(q^{2n+*}\Big)^{-s},$$
which is easily seen to be convergent if ${\rm Re}(s)>1$ and rational as well.
In this way, we complete the proof of the Theorem 3.
 
In fact,  a functional equation also holds for these new zeta functions. 
But we will
leave it to the reader. Clearly, from (the proof of) this theorem, we may 
conclude
that  zeta functions  associated to non-stable bundles, and hence the so-called
total zeta function associated to all bundles, both semi-stable and
non-stable, are interesting but not as canonical as that for semi-stable 
bundles.
Moreover, note that total zeta functions are closed related to  automorphic
$L$-functions. Hence, it would be of great interests if one  does  similar
decompositions (as what we have done here) for them.

We end this paper by pointing out that this work is a continuation
of our small note [We1], in which a different kind of restricted zeta 
function is
formulated, and  that all these works are motivated by our
new  non-abelian zeta functions for number fields [We2].
\vskip 0.45cm
\centerline {\bf REFERENCES}
\vskip 0.30cm
\item{[A]} E. Artin, {\it Collected Papers}, edited by S. Lang and J. Tata,
Addison-Wesley (1965)
\vskip 0.30cm
\item{[AB]} M.F. Atiyah \& R. Bott, The Yang-Mills equations over Riemann 
surfaces,
Phil. Trans. R. Soc. Lond. A. {\bf 308}, 523-615 (1982)
\vskip 0.30cm
\item{[DR]} U.V. Desale \& S. Ramanan, Poincar\'e polynomials of the variety of
stable bundles, Maeh. Ann {\bf 216}, 233-244 (1975)
\vskip 0.30cm
\item{[HN]} G. Harder \& M.S. Narasimhan, On the cohomology groups of 
moduli spaces
of vector bundles over curves, Math Ann. {\bf 212}, (1975) 215-248
\vskip 0.30cm
\item{[Mu]} D. Mumford, {\it Geometric Invariant Theory}, Springer-Verlag, 
Berlin
(1965)
\vskip 0.30cm
\item{[Se]} J.P. Serre, {\it Trees}, Springer-Verlag, Berlin (1980)
\vskip 0.30cm
\item{[W]} A. Weil, {\it Sur les courbes alg\'ebriques et les vari\'et\'es qui
s'en d\'eduisent}, Herman, Paris (1948)
\vskip 0.30cm
\item{[We1]} L. Weng, A result on zeta functions for curves over finite fields,
preprint, Nagoya University, 2000
\vskip 0.30cm
\item{[We2]} L. Weng,  Riemann-Roch Theorem, Stability, and New Zeta
Functions for Number Fields, preprint, Nagoya, 2000
\end